%% file: root.tex
\newif\ifcomments

\newif\ifproofs
\proofstrue

\documentclass[letterpaper, 10 pt, conference]{ieeeconf}  

\IEEEoverridecommandlockouts                              
\usepackage[english]{babel}
\usepackage{verbatim}
\usepackage{graphicx}
\usepackage{amsfonts}
\usepackage{amsmath}
\usepackage{amssymb}
\usepackage{paralist}
\usepackage{wrapfig}
\usepackage{xcolor}
\usepackage{flushend}

\usepackage{subcaption}

\usepackage{siunitx}

\definecolor{agro}{HTML}{228B9D}

\usepackage[capitalise, nameinlink]{cleveref}

\usepackage{pgfplots}

\graphicspath{{imgs/}{../experiments/forgetfull_bouncing_ball/}}
\newcommand{\defaultImageWidth}{0.3}
\newcommand{\defaultHalfImageWidth}{0.3}

\input{commands}

\newcommand{\nodelayorbitfig}{\cite[Figure~6.1]{Gomes2019d}}
\newcommand{\delayorbitfig}{\cite[Figure~6.2]{Gomes2019d}}
\newcommand{\figbimodalsystemdelayedswitching}{\cite[Figure~6.4]{Gomes2019d}}
\newcommand{\fighplot}{\cite[Figure~6.14]{Gomes2019d}}

\bstctlcite{rmurl}

\newcommand{\ceqref}[1]{\eqref{#1}}
\newcommand{\ceqreff}[2]{\eqref{#1}, \eqref{#2}}

\crefname{problem}{problem}{problems}

\title{\LARGE \bf
Stability of Planar Switched Systems\\ under Delayed Event Detection
}

\author{Beno\^it Legat$^{1}$, Cl\'audio Gomes$^{2}$, Paschalis Karalis$^{3}$, \\Rapha\"el M. Jungers$^{4}$, Eva M. Navarro-L\'opez$^{5}$, and Hans Vangheluwe$^{6}$%
	\thanks{$^{1}$B. L. is a F.R.S.-FNRS Research Fellow, at UCLouvain.}%
	\thanks{$^{2}$C. G. is a FWO Research Fellow, at the University of Antwerp, supported by the Research Foundation - Flanders (File Number 1S06316N).}%
	\thanks{$^{3}$P. K. is a research engineer at Cadence Design Systems.}%
	\thanks{$^{4}$R. M. Jungers is a F.R.S.-FNRS honorary Research Associate, at UCLouvain. He is supported by the Walloon Region and the Innoviris Foundation.}%
	\thanks{$^{5}$E. M. N. is a Honorary Senior Research Fellow at the University of Manchester.}%
	\thanks{$^{6}$H. V. is a Professor at the University of Antwerp. His work is partially supported by Flanders Make vzw, the strategic research centre for the manufacturing industry.}%
}

\begin{document}

\maketitle
\thispagestyle{empty}
\pagestyle{empty}

\input{abstract}

\input{introduction}

\input{example}

\input{background}

\input{orbit}

\input{application}

\input{conclusion}

\bibliographystyle{IEEEtran}
\bibliography{rmurl,bibliography_claudio,bibliography_benoit}

\ifproofs %
  \appendix
\else %
\fi%

\ifproofs %
  \input{proofs}
\else %
\fi%

\end{document}

%% file: commands.tex
\newcommand{\claudio}[1]{%
\ifcomments %
{\scriptsize
\textbf{\textcolor{red}{Claudio: #1}}
} %
\else %
\fi%
}%

\newcommand{\benoit}[1]{%
\ifcomments %
{\scriptsize
\textbf{\textcolor{blue}{Benoit: #1}}
} %
\else %
\fi%
}%

\newcommand{\paschalis}[1]{%
\ifcomments %
{\scriptsize
\textbf{\textcolor{magenta}{Paschalis: #1}}
} %
\else %
\fi%
}%

\newcommand{\eva}[1]{%
\ifcomments %
{ \scriptsize
\textbf{\textcolor{purple}{Eva: #1}}
} %
\else %
\fi%
}%

\newtheorem{assumption}{Assumption}{}

\crefname{assumption}{assumption}{assumptions}

\DeclareMathOperator*{\interior}{int}

\newcommand{\HS}[0]{\ensuremath{\mathcal{S}}}

\newcommand{\SD}[0]{\ensuremath{\HS}} 

\DeclareMathOperator{\sign}{sign}
\newcommand{\swdl}{delay in the transition} 

\newcommand{\MSD}{MSD}
\newcommand{\MSDdef}{Maximum Stable Delay (\MSD{})}
\newcommand{\trdl}{transition-delayed}
\newcommand{\trdlsys}{transition-delayed system}

\newcommand{\setvar}{U}
\newcommand{\Hsys}{Hybrid automaton}
\newcommand{\hsys}{hybrid automaton}
\newcommand{\hsyss}{hybrid automata}
\DeclareMathOperator{\Init}{Init}
\DeclareMathOperator{\Dom}{Dom}
\DeclareMathOperator{\Guard}{Guard}

\newcommand{\execvar}{\mathcal{T}}
\newcommand{\timeslice}[3]{{#1}{[#2, #3]}}
\newcommand{\timeexec}[2]{\timeslice{\execvar}{#1}{#2}}
\newcommand{\stateproj}[1]{{#1}|_{\mathcal{X}}}
\newcommand{\stateexec}{\stateproj{\execvar}}
\newcommand{\statetimeexec}[2]{\stateproj{\timeslice{\execvar}{#1}{#2}}}
\newcommand{\hmapexec}[2]{\execvar\hpcmap{#1}{#2}}
\newcommand{\statehmapexec}[2]{\stateproj{\execvar\hpcmap{#1}{#2}}}
\newcommand{\joinexec}[2]{{#1} \uplus {#2}}
\newcommand{\anhexec}{a hybrid trajectory}
\newcommand{\anexec}{a trajectory}
\newcommand{\hexec}{hybrid trajectory}
\newcommand{\exec}{trajectory}
\newcommand{\execs}{trajectories}
\newcommand{\PCMap}{Poincar\'e Map}
\newcommand{\PCRay}{Poincar\'e curve}

\newcommand{\TimeModeTwoSp}{T_2}

\newcommand{\PointCarreLine}{S_p}
\newcommand{\GoodPointCarreLine}{\tilde{S}_p}
\newcommand{\pcfun}{s}
\newcommand{\SwitchingSurface}{G}
\newcommand{\GoodSwitchingSurface}{\tilde{\SwitchingSurface}}
\newcommand{\swfun}{g}
\newcommand{\GoodD}{D}
\newcommand{\rtime}[3]{\tau_{#2}^{#1}(#3)}
\newcommand{\fullpcmap}[3]{\mathcal{P}_{#2}^{#1}(#3)}
\newcommand{\pcmap}[3]{P_{#2}^{#1}(#3)}
\newcommand{\execpcmap}[3]{\execvar\pcmap{f_{#1}}{#2}{#3}}
\newcommand{\hpcmap}[2]{\mathcal{Q}(#2, #1)}
\newcommand{\Floww}{\phi}

\newcommand{\Flowf}[3]{\phi_{#1}(#2, #3)}
\newcommand{\execFlow}[3]{\execvar\Flowf{f_{#1}}{#2}{#3}}
\newcommand{\Flowi}[3]{\Flowf{f_{#1}}{#2}{#3}}

\newtheorem{definition}{Definition}{}
\newtheorem{example}{Example}{}
\newtheorem{lemma}{Lemma}{}
\newtheorem{remark}{Remark}{}
\newtheorem{corollary}{Corollary}{}
\newtheorem{proposition}{Proposition}{}
\newtheorem{problem}{Problem}{}
\newtheorem{theorem}{Theorem}{}
{}
{}

\newcommand{\opt}[1]{\ensuremath{#1^\star}}
\newcommand{\closure}[1]{\ensuremath{\overline{#1}}}
\newcommand{\lcinter}[1]{\ensuremath{ \left[ #1 \right) }}
\newcommand{\brackets}[1]{\ensuremath{ \left[ #1 \right] }}

\newcommand{\set}[1]{\ensuremath{ \left\{ #1 \right\}}}

\newcommand{\dert}[1]{\ensuremath{ \dot{#1} }}

\newcommand{\setreal}[0]{\ensuremath{\mathbb{R}}}

\newcommand{\norm}[1]{\left\lVert#1\right\rVert}
\newcommand{\enorm}[1]{\norm{#1 - x^*}_2} 

\newcommand{\vectorOne}[1]{\brackets{%
\begin{matrix}
  #1
 \end{matrix}%
}}

\newenvironment{aligneq*}%
{
\begin{equation*}
\begin{aligned}
}{
\end{aligned}
\end{equation*}
}

\newenvironment{aligneq}%
{
\begin{equation}
\begin{aligned}
}{
\end{aligned}
\end{equation}
}


%% file: abstract.tex
\begin{abstract}


In this paper, we analyse the impact of delayed event detection
on the stability of a 2-mode planar hybrid automata.
We consider hybrid automata with a unique equilibrium point for all the modes,
and we find the maximum delay that preserves stability of that equilibrium point.
We also show for the class of hybrid automata treated that the instability of the equilibrium point for the equivalent hybrid automaton with delay in the transitions is equivalent to the existence of a closed orbit in the hybrid state space, a result that is inspired by the Joint Spectral Radius theorem.
This leads to an algorithm for computing the maximum stable delay \emph{exactly}.
Other potential applications of our technique include co-simulation, networked control systems and delayed controlled switching with a state feedback control.
\end{abstract}

%% file: introduction.tex
\section{Introduction}
\label{sec:introduction}

In order to make the computation of the behaviour of hybrid systems possible, one needs to not only approximate the continuous dynamics, but also to \emph{accurately} and \emph{efficiently} detect when to compute the discrete-event dynamics.
In order to accurately do so, transition (or event) detection and location schemes are employed.
The efficiency requirement is satisfied by setting the appropriate parameters of the transition location scheme. Wrong tolerance values can lead to unnecessary computations and/or inaccurate results (see, e.g., \cref{sec:motivating_example}).

\textbf{Contributions.}
In this paper, we focus on the numerical stability of hybrid system simulation.
We formalize the problem of finding the \MSDdef{}, that is, the maximum delay a hybrid system can admit in a transition, so that its solutions remain asymptotically stable (\cref{sec:problem_form}).

Inspired by the Joint Spectral Radius theorem \cite{Jungers2009}, we show how this problem can be solved for a restricted class of planar hybrid systems, by reducing it to a problem of finding a non-trivial closed orbit in the \trdlsys{} (\cref{sec:orbit}).
As an example application of these theoretical results, we provide an algorithm to find such a closed orbit in a non-linear hybrid system (\cref{sec:result}).

The \MSD{} can then be used to set the appropriate parameterization of transition detection and location schemes, and/or the size of the simulation step.
Other potential applications include:
\begin{inparaenum}[(i)]
\item networked control systems (see, e.g., \cite{Bauer2012}), where the component responsible for deciding the mode of the plant may be reacting to a delayed signal; and
\item real-time simulation of hybrid systems, where state transition location can be disabled, or relaxed, to increase performance, as in animation of colliding multi-body systems~\cite{Conti2016}.
\end{inparaenum}

An extended version of this work is available in \cite[Chapter 6]{Gomes2019d}.

\textbf{Related Work.}
In the domain of network control, \cite{Jungers2012} focuses on the state-feedback stabilization of LTI systems that can be remotely controlled over a multi-hop network. In contrast, we focus on non-linear planar systems with a switching surface.
We follow an approach that is similar to \cite{Gomes2017c,Bauer2012,Giannakopoulos2001}.
Prior work \cite{Gomes2017c} explores how to use Lyapunov functions to approximate the \MSD{} in the particular case of affine systems, without resorting to finding closed orbits.
However, no formal proof is provided regarding the correctness of the procedure.
In \cite{Bauer2012}, the delayed system is formulated as a hybrid system. Then, the stability of which is proven if a Lyapunov function can be found, using the sum of squares \cite{Papachristodoulou2005} approach, which can only approximate the \MSD{}.
The work in \cite{Giannakopoulos2001} provides a comprehensive study of the dynamics of planar systems with the form $\dert{x} = Ax + \sign(w^T x)v$, with $A$ a constant matrix.
Similarly to the current manuscript, Poincaré maps are used to analyze the stability of periodic motions (\cite[Section 5.2]{Giannakopoulos2001}).
More recently, \cite{proskurnikovdoes} studied the maximum sampling time $\tau$ such that the asymptotic stability is preserved under the discretization of a hybrid nonlinear system
given that the length of the (not necessarily equal) sampling intervals does not exceed $\tau$.
While the problem setting is similar to our work, their aim is to prove the existence of a maximal sampling time $\tau$, not to approximate it.
We refer the reader to \cite{Arapostathis2007} for similar work restricted to linear.

%% file: example.tex
\section{Motivating Example: Relaxed Bouncing Ball}
\label{sec:motivating_example}

Consider a bouncing ball modeled with two modes: free-fall, and contact. The detailed equations for each mode are given in \cite[Section 6.2]{Gomes2019d}.

The ball changes from free fall mode to contact mode when it comes in contact with the floor.
Formally, that is when $g(x) = \vectorOne{1 & 0} x - r \leq 0$, with $r$ denoting the radius of the ball.

\nodelayorbitfig{} shows an example \exec{}, which converges to an equilibrium in the contact mode. 
In a simulation, the moment that the ball changes from the free-fall mode to the contact mode is approximated using a state transition location technique.
\delayorbitfig{} shows a simulation where the transition from free-fall to contact mode is not correctly detected, resulting in a delay of $0.002$s.
The simulation shows the ball reaching the same compression that it had on its initial state.
Since the total energy of the physical bouncing ball
is dissipated via air friction and impact damping, the bouncing ball simulation should come to a rest.
Instead, \delayorbitfig{} shows that the simulation of the bouncing ball with transition delays of up to $0.002$s may fail to come to a rest.

The \textbf{research question} follows:
for a given hybrid system whose trajectories eventually tend to an equilibrium, what is the largest \swdl{} the simulations can tolerate, such that that property is preserved?
The next section formulates this problem.

%% file: background.tex
\section{Problem Formulation}
\label{sec:problem_form}

\input{background_hybrid_systems}
We say that an equilibrium of a \hsys{} is \emph{Globally asymptotically stable (GAS)} if all \execs{} of the system
remain close, and
converge, to the equilibrium (see, e.g., \cite{Goebel2012}).

As motivated in \cref{sec:motivating_example}, we need to analyze the behavior of a \hsys{} under \swdl{}.
It turns out that we can construct a \hsys{} for which the stability is equivalent to the stability under \swdl{} of the original system.
We present this reduction in the following restricted class of \hsyss{} depicted in \cref{eq:bimodal_system}.

\begin{definition}[Bi-modal Hybrid Automata]\label{def:bms}
  A \emph{Bi-modal Hybrid Automaton} (BMHA) system is a \hsys{} as defined in \cref{def:ha} such that
  there exists a continuously differentiable function $\swfun$ such that
  $Q = \{1, 2\}$,
  $\Dom(1) = \{\, x \mid \swfun(x) \le 0 \,\}$,
  $\Dom(2) = \{\, x \mid \swfun(x) \ge 0 \,\}$,
  $E = \{(1, 2), (2, 1)\}$,
  $\Guard((1, 2)) = \{\, x \mid \swfun(x) \ge 0 \,\}$, $\Guard((2, 1)) = \{\, x \mid \swfun(x) \le 0 \,\}$
  and $R((1, 2), x) = R((2, 1), x) = \{x\}$ for all $x \in \setreal^n$.
  This is illustrated in \cref{eq:bimodal_system}.
\end{definition}

\begin{definition}
  \label{def:SDH}
  Given a BMHA $\HS$, as defined in \cref{def:bms}, and a delay $H > 0$, we define its \trdl{} $\SD_H$ counterpart as in \figbimodalsystemdelayedswitching{}.
\end{definition}


$\SD_H$ is non-deterministic, having \execs{} where a transition can be delayed by a maximum of $H$ units of time.
Moreover, $\SD_H$ has the same equilibrium point as the original BMHA system.
And the equilibrium point of $\SD_H$ may not be GAS, even though the original system $\HS$ is GAS, as exemplified in \delayorbitfig{}.

The following results follow from \cref{def:bms,def:SDH}, and lead to our problem formulation.

\begin{proposition}\label{prop:largerH_covers_everything}
  If $\execvar$ is \anexec{} of $\SD_{H'}$, then it is also \anexec{} of $\SD_{H}$, for any $H' \leq H$.
\end{proposition}

\begin{corollary}\label{stability_smaller_H_contrapos}
For any $H' \leq H$, if the equilibrium $x^* \in \setreal^n$ of $\SD_H$ is GAS, then the equilibrium $x^*$ of $\SD_{H'}$ is also GAS.
\end{corollary}

The previous results naturally prompt the following problem.
\begin{problem}\label{eq:general_problem}
  Given a BMHA $\HS$, as defined in \cref{def:bms}, with a GAS equilibrium $x^*$, find
\begin{aligneq}
\sigma(\SD) = \sup_{H \geq 0} H \text{ s.t. } x^* \text{ is a GAS equilibrium of } \SD_H
\end{aligneq}
\end{problem}
We denote the \emph{switching surface} as
\begin{equation*}
\SwitchingSurface = \{\, x \in \mathbb{R}^2 \mid g(x)=0 \,\}.
\end{equation*}
To avoid pathological cases, we assume that the single equilibrium is not in the switching surface.
Furthermore, without loss of generality, we consider that the equilibrium is in the interior of
$\Dom(1)$ and assume that it is GAS for the mode 1.

We call $\sigma(\SD)$ the \MSDdef{}.
Intuitively, $\sigma(\SD)$ gives us the maximal delay for which the GAS property $x^*$ is preserved under delayed switching:
by \cref{stability_smaller_H_contrapos}, $x^*$ of $\SD_H$ is GAS for all $H < \sigma(\SD)$.
Knowing $\sigma(\SD)$ allows one to define the time step of a simulation algorithm, or the tolerance of a state transition location scheme.

\Cref{eq:general_problem} can be extended to a general \hsys{} (\cref{def:ha}) because the construction of $\SD_H$ can be generalized.
We focus on BMHA because our subsequent results cannot be straightforwardly generalized. This is left as future work.

\Cref{eq:general_problem} can be reduced to the stability of an equilibrium of a \hsys{} of the form given in \figbimodalsystemdelayedswitching{},
and thus we can leverage classical tools developed for the stability of \hsyss{}.
For example, \cite[Section 6.5.1]{Gomes2019d} demonstrates the use of SpaceEx to prove the stability of $\SD_H$, for a given $H$.
However, deciding the stability of an equilibrium of a \hsys{} is not possible in general \cite{blondel1999complexity}, 
hence this does not provide an \emph{efficient} approach and it only \emph{approximates} the $\sigma(\SD)$ of \cref{eq:general_problem}.

As we show in the following section,
under mild assumptions, the stability of the subclass of planar \hsyss{} of the form given in \figbimodalsystemdelayedswitching{}
is equivalent to the existence of a closed orbit (\cref{def:orbit}).
As shown in \cref{sec:result}, this theoretical result leads to \emph{practical} algorithms
to compute $\sigma(\SD)$ with \emph{arbitrary accuracy}.

%% file: background_hybrid_systems.tex
We adopt the usual definition of Hybrid Automaton and the notation in \cite[Definition 26]{Gomes2019d}.

\begin{definition}[\Hsys{} \cite{johansson1999regularization}]
\label{def:ha}
A \hsys{} $\HS$ is a collection
$$\HS = (Q, E, \mathcal{X}, \Dom, \mathcal{F}, 
  \Guard, R)$$
where:
\begin{inparaitem}[\textbullet]
\item $Q$ is a finite set of modes;
\item $E \subseteq Q \times Q $ is a finite set of edges called transitions;
\item $\mathcal{X} \subseteq \setreal^n$ is the continuous state space, for some natural $n$;
\item $\Dom : Q \to 2^{\mathcal{X}}$ is the mode domain;
\item $\mathcal{F} = \set{ f_{q} (x) : q \in Q }$ is a collection of time-invariant vector fields such that each $f_{q} (x)$
is Lipschitz continuous on $\mathit{Dom}(q)$;
\item $\Guard : E \to 2^\mathcal{X}$ defines a guard set for each transition. 
\item $R : E \times \mathcal{X} \to 2^\mathcal{X}$ specifies how the continuous state is reset at each transition.
\end{inparaitem}
\end{definition}

The Lipschitz continuity assumption in \cref{def:ha} ensures the existence and continuity of a solution inside a mode as shown by \cite[Lemma 1]{Gomes2019d}.

\begin{definition}[Flow]
  Given a dynamical system $\dot{x}(t) = f(x(t))$,
  the \emph{flow of the vector field} $f$,
  denoted $\Flowf{f}{x}{t}$,
  is the state of the system after a time $t$ if the initial state is $x$.
\end{definition}

\begin{figure}[htp]
\begin{center}
  \includegraphics[width=0.4\textwidth]{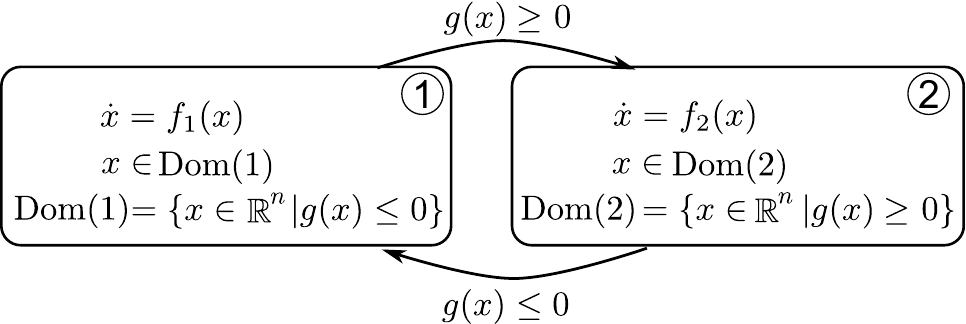}
  \caption{A Bi-modal Hybrid Automata in the notation used by \cite{Johansson1999}.}
  \label{eq:bimodal_system}
\end{center}
\end{figure}

\begin{example}\label{ex:bouncing_ball}
  \cref{eq:bimodal_system} shows the graphical representation of the hybrid automaton for the bouncing ball described in \cref{sec:motivating_example}.
\end{example}

\begin{definition}[{Hybrid \exec{}.  \cite[Definition~2.3]{Goebel2012}}]
\label{def:hybrid_time}
  A hybrid \exec{} of a system $\HS$ as defined in \cref{def:ha} is a
  subset $\execvar \subseteq Q \times \setreal \times \mathcal{X}$ with
	\footnotesize
  \[
	\execvar = \bigcup_{j=0}^{J-1} \{q_j\} \times \Big\{\, (t, x) \mid t_j \le t \le t_{j+1}, x = \Flowf{f_{q_j}}{x_j}{t-t_j} \in \Dom(q_j) \,\Big\}
 \]
	\normalsize
  for some finite sequence of times $t_0 < t_1 < \cdots < t_J$,
  modes $q_j \in Q$ and state vector $x_j \in \Dom(q_j)$ such that
  $(q_{j-1}, q_j) \in E$,
  $y_j \in \Guard((q_{j-1}, q_j))$
  and $x_j \in R((q_{j-1}, q_j),y_j)$ for $j = 1, \ldots, J$
  where $y_j = \Flowf{f_{q_{j-1}}}{x_{j-1}}{t_j-t_{j-1}}$.
\end{definition}

Note that \cref{def:hybrid_time} excludes chattering and Zeno behaviour, which is assumed inexistent later by \cref{ass:perppc} but considered in \cite[Definition~2.3]{Goebel2012}.
The later assumption also allows restricts the non-determinism of hybrid automata while still permitting non trivial behaviour.

We define the following notations for \execs{}:
\begin{align*}
  \timeexec{\underline{t}}{\overline{t}} & = \{\, (q, t, x) \in \execvar \mid \underline{t} \le t \le \overline{t} \,\},\\
  \stateexec & = \{\, x \mid \exists q, t: (q, t, x) \in \execvar \,\}\\
  t_{\text{min}}(\execvar) & = \inf_{(q, t, x) \in \execvar} t \qquad \qquad t_{\text{max}}(\execvar) = \sup_{(q, t, x) \in \execvar} t\\
  \joinexec{\execvar_1}{\execvar_2} & = \execvar_1 \cup \{\, (q, t + t_{\text{max}}(\execvar_1) - t_{\text{min}}(\execvar_2), x) \mid (q, t, x) \in \execvar_2 \,\}.
\end{align*}
Note that if $\execvar$ is a hybrid \exec{} then so is $\timeexec{\underline{t}}{\overline{t}}$,
which allows us to define the combined notation $\statetimeexec{\underline{t}}{\overline{t}}$.
However, even if $\execvar_1$ and $\execvar_2$ are hybrid \execs{},
$\joinexec{\execvar_1}{\execvar_2}$ may not be a hybrid \exec{} if
for instance it induces a transition that does not respect a guard set.

%% file: orbit.tex
\section{Orbit and stability}
\label{sec:orbit}

Our approach is similar to the stabilization of a \emph{constrained linear switched system} \cite{Dai2012,Philippe2016}.
In the theory of constrained switched systems, the asymptotic stability, can be proven by showing that there are no closed orbits in the system (see \cite[Theorem A]{Dai2012}, \cite{Lagarias1995}, and example applications in \cite{Gomes2018c}).
We show that there is an equivalence between the absence of cyclic behaviors in the \trdl{} counterpart (\cref{def:SDH}), and its asymptotic stability.
We start by defining the notion of closed orbit, to formalize cyclic \execs{}.

\begin{definition}[Closed orbit]
  \label{def:orbit}
  Given \anhexec{} $\execvar$ as given in \cref{def:hybrid_time}, 
  we say that it is a closed orbit if there exists
  $$(q_0, t_0, x_0), (\tilde{q}, \tilde{t}, \tilde{x}), (\bar{q}, \bar{t}, \bar{x}) \in \execvar, $$
  such that
  $\bar{q} = q_0$,
  $t_0 < \tilde{t} < \bar{t}$
  and $\tilde{x} \neq \bar{x} = x_0$.
\end{definition}
The condition $\tilde{x} \neq x$ excludes the equilibrium from being a trivial closed orbit.

\begin{proposition}\label{claim:orbit}
  If there is a closed orbit in $\SD_H$, then the equilibrium $x^*$ of $\SD_H$ is not GAS.
\end{proposition}

The above proposition allows, for a given $H$, to prove the instability of $x^*$ of $\SD_H$ by just finding a closed orbit.
By trying to find the smallest $H$ for which there is a closed orbit, one can approximate from above the solution $\sigma(\SD)$ to \cref{eq:general_problem}.
However, approximation from above may not be sufficient.
To see why, imagine that the smallest $H$ for which there exists a closed orbit in $\SD_H$ is found.
Then, it might still be the case that there exists a $H' < H$ such that $x^*$ in $\SD_{H'}$ is not GAS.
Our contribution is to show that such $H'$ cannot exist for a class of planar BMHA systems
satisfying \cref{ass:closedinj,ass:perppc,ass:perpsw,ass:perpsw1}, formalized below.

We use classical results of analysis of planar non-linear systems such as continuity and monotonicity of a \emph{\PCMap{}}.
Such classical notions need to be carefully used in this setting since they are usually developed
for non-hybrid systems.
We start by discussing the monotonicity.
We will see in \cref{coro:decreasing} that the \PCMap{} is monotonous
for the \hsys{} $\SD$.
However, for a non-deterministic \hsys{} such as $\SD_H$, monotonicity of the \PCMap{} is not guaranteed
as two \execs{} can intersect as long as they are in different modes.
Nevertheless, the topological argument commonly used to prove monotonicity (see e.g. \cite[Section~10.4]{Hirsch2012})
can still be used to obtain the result given in \cref{lem:intersect} for planar BMHA systems.

We start by defining the following connectedness concepts from topology;
see e.g. \cite[Chapter~3]{Munkres2000}.

\begin{definition}[Path]
  \label{def:path}
  Given a set $X$ and two points $x, y \in X$, a \emph{path} in $X$
  from $x$ to $y$ is the image of a continuous map $f: [0, 1] \to X$
  such that $f(0) = x$ and $f(1) = y$.
  We denote the union of all paths in $X$ from $x$ to $y$ as $[x, y]_X$.
\end{definition}

\begin{definition}[Path components of $X$]
  \label{def:path_component}
  We define an equivalence relation on the set $X$ by defining
  $x \sim y$ if there is a path in $X$ from $x$ to $y$.
  The equivalence classes are called the \emph{path components} of $X$.
\end{definition}

Given a set $\setvar$, we denote its \emph{closure} by $\closure{\setvar}$
and its \emph{interior} by $\interior(\setvar)$.

Given the vector field $f_1$ of the first mode of a planar BMHA $\SD$ defined in
\cref{def:bms},
we assume the existence of a continuously differentiable function $\pcfun : \setreal^2 \to \setreal$
such that the curve
\begin{equation}
  \label{eq:pc}
  \PointCarreLine = \{\, x \mid \pcfun(x) = 0, \nabla \pcfun(x) \cdot f_1(x) > 0 \,\}
\end{equation}
satisfies $\PointCarreLine \cap \SwitchingSurface = \emptyset{}$, where $\SwitchingSurface$ is the switching surface.
We refer to the closed section $\PointCarreLine$ as \PCRay{},
its definition is similar to the notion of \emph{local sections}; see \cite[Section~10.2]{Hirsch2012}.
\begin{assumption}
  \label{ass:closedinj}
  The \PCRay{} $\PointCarreLine$ has a single path component,
  $\closure{\PointCarreLine} = \PointCarreLine \cup \{x^*\}$
  and the restriction of the Euclidean norm to $\PointCarreLine$ is injective.
\end{assumption}
We use the following notations for the \emph{return time} and \emph{\PCMap{}}, given a vector field $f$ and a set $\setvar$:
\begin{align*}
  \rtime{f}{\setvar}{x} & = \inf \{\, t \mid t > 0, \Flowf{f}{x}{t} \in \setvar \,\}, \\
  \pcmap{f}{\setvar}{x} & = \Flowf{f}{x}{\rtime{f}{\setvar}{x}}.
\end{align*}
Note that $\pcmap{f}{\setvar}{x}$ is not defined when $\rtime{f}{\setvar}{x}$ is infinite.

Given a \exec{} $\execvar$, a set $\setvar$ and a time $t_0 \in \setreal$,
we define the following notation
\begin{align*}
  \rtime{\execvar}{\setvar}{t_0} & = \inf \{\, t \mid (q, t, x) \in \execvar, t > t_0, x \in \setvar \,\}, \\
  \fullpcmap{\execvar}{\setvar}{t_0} & = \{\,(q, t, x) \in \execvar \mid t = \rtime{\execvar}{\setvar}{t_0} \,\}.
\end{align*}



The following lemma shows that, whenever a \exec{} $\execvar$ of $\SD_H$ intersects with itself, it must be under a different mode.
\cite[Figure 6.5]{Gomes2019d} gives representative examples of this situation.

\begin{lemma}
  \label{lem:intersect}
  Given a planar \hsys{} $\SD_H$ as given in \cref{def:SDH}, with a \PCRay{} $\PointCarreLine$ as given in \ceqref{eq:pc}
  and a \exec{} $\execvar$ of $\SD_H$, 
  consider 
  $(q_1, t_1, x_1) \in \execvar$
  with $x_1 \in \PointCarreLine$
  and $(q_2, t_2, x_2) \in \fullpcmap{\execvar}{\PointCarreLine}{t_1}$
  with $\enorm{x_2} \geq \enorm{x_1}$.
  Under \cref{ass:closedinj}, there exists a path component $\setvar$ of $\mathcal{X} \setminus (\statetimeexec{t_1}{t_2} \cup [x_2, x_1]_{\PointCarreLine})$
  such that:
  \begin{compactenum}[(1)]
    \item $\{\, x \in \PointCarreLine \mid \enorm{x} \leq \enorm{x_2} \,\} \subseteq \closure{\setvar}$; and
    \item there exists $\epsilon > 0$ such that $\statetimeexec{t_2}{t_2 + \epsilon} \cap U = \emptyset$.
  \end{compactenum}
  Moreover, if $\statetimeexec{t_2}{\infty} \cap \setvar \neq \emptyset$ then the intersection
  $(q_3, t_3, x_3) \in \fullpcmap{\execvar}{\closure{\setvar}}{t_2}$ is such that
  there exists $(q_4, t_4, x_3) \in \timeexec{t_1}{t_2}$ with $q_3 \neq q_4$.
\end{lemma}

\paschalis{It appears to me that it is enough to say that none of the modes has cycles of its own -- ie, that none of the subsystems (say x = fq(x)) has periodic trajectories. This is a mild assumption, i think -- it is also easy to describe without paths and path connectedness). With it, Lemma 2 follows trivially, i think.}
\benoit{You need to prove that if you have not crossed the switching surface between two intersection of the poincarre map, you won't cross it again. If you don't prove that, you cannot prove that you stay ini mode 1 so saying that mode 1 has no orbit is not enough. The reason you won't cross it again is because you have created a path component in which you are now stuck inside and this path component does not intersect with the switching surface. Therefore proving Lemma 2 with Lemma 1 seems natural.}

\paschalis{In general, we have defined paths as just continuous functions -- it is not clear why we care if there is a path component (a subset where we can connect stuff with paths).}
\benoit{It's just an elegant way (that is handy in this case since we have already defined path) to define the intuitive notion of ``a set in one part''. It's like connected compontents for graphs. It's a bit like convex sets buts it's more general, you can have holes in path componenents.}

\paschalis{The assertion seems to be that there is a connected subset of Sp, for which the next intersection with Sp under SH is always closer to x*, it it never happens in 0 time -- is that correct?}
\benoit{Yes exactly, it's proving this using a topologically argument relying on the planarity of the system}


As a corollary, the monotonicity of the \PCMap{} holds for the mode 1, as stated in \cref{coro:decreasing}.
However, when applying the same result to the delayed \hsys{}, care must be taken to consider only trajectories that remain in mode 1, as these are not affected by the delayed switching.
If a trajectory of the delayed \hsys{} remains in the domain of mode 1 between two consecutive
intersections with the \PCMap{},
\cref{coro:decreasing} shows that the second intersection has a smaller Euclidean norm.
In fact, as shown by \cref{lem:invariant_mode_1},
this implies that the trajectory will remain in mode 1 until its end.

\begin{corollary}
  \label{coro:decreasing}
  Consider a planar \hsys{} $\SD_H$ as given in \cref{def:SDH}
  and a \PCRay{} $\PointCarreLine$ as defined in \ceqref{eq:pc}.
  For all $x \in \PointCarreLine$,
  if for any $t$, $\Flowf{f}{x}{t} \in \Dom(1)$, then $\enorm{\pcmap{f_1}{\PointCarreLine}{x}} < \enorm{x}$.
\end{corollary}

\begin{lemma}
  \label{lem:invariant_mode_1}
  Given a delay $H$ and a planar \hsys{} $\SD_H$ as given in \cref{def:SDH}, with a
  \PCRay{} $\PointCarreLine$ as defined in \ceqref{eq:pc},
  and a \exec{} $\execvar$ of $\SD_H$, 
  consider 
  $(1, t_1, x_1) \in \execvar$
  with $x_1 \in \PointCarreLine$
  and $(1, t_2, x_2) \in \fullpcmap{\execvar}{\PointCarreLine}{t_1}$.
  Under \cref{ass:closedinj},
  if $\statetimeexec{t_1}{t_2} \subseteq \Dom(1)$, then
  $\statetimeexec{t_2}{\infty} \subseteq \Dom(1)$.
\end{lemma}

In view of \cref{coro:decreasing} and \cref{lem:invariant_mode_1},
if \anexec{} does not intersect the switching surface between two intersections of
the \PCRay{} (hence staying in mode 1) then it will repeat this behavior
indefinitely. 
Therefore, as our aim is to find a closed orbit, we restrict our attention to
points of the \PCRay{} from which \anexec{} of $\SD$ intersects the switching
surface.
We denote the set of such points, illustrated in \cite[Figure 6.6]{Gomes2019d}, as follows:
\begin{multline}\label{eq:goo_pc_line}
  \GoodPointCarreLine = \{\, x \in \PointCarreLine \mid \rtime{f_1}{\SwitchingSurface}{x}, \text{ is finite},\\
  \nabla g(\pcmap{f_1}{\SwitchingSurface}{x}) \cdot f_1(\pcmap{f_1}{\SwitchingSurface}{x}) > 0 \,\}.
\end{multline}
Note that, by continuity of the functions $f_1$ and $g$, the value $g(\pcmap{f_1}{\SwitchingSurface}{x}) \cdot f_1(\pcmap{f_1}{\SwitchingSurface}{x})$ in \eqref{eq:goo_pc_line} cannot be negative if $x \in \Dom(1)$,
hence the condition $> 0$ could be replaced by $\neq 0$.
The condition simply excludes points $x$ such that $f_1(\pcmap{f_1}{\SwitchingSurface}{x})$
is tangent to the switching surface.

Since the equilibrium $x^*$ is GAS for the mode 1,
\paschalis{This has to be an assumption -- even if we can show that it follows by the sequence:
  assumption 1 --> lemma 2 --> corollary 2, saying that x* is GAS for f1 is much clearer.}
by \cref{coro:decreasing} and \cite[Lemma 1]{Gomes2019d},
any trajectory of $\SD_H$ (independently of $H$) starting at a point in
$\PointCarreLine \setminus \GoodPointCarreLine$ will converge to $x^*$.
\cref{lem:U} provides a converse result, that is,
if $x^*$ is not GAS and $\SD_H$ does not admit any closed orbit,
then subsequent intersections with $\PointCarreLine$ will be farther from the origin.
However, extra care must be taken for where to place the \PCRay{}, to make sure all its intersections with $\PointCarreLine$ happen in mode 1.
Moreover, we need to assume that the dynamics of mode 2, under the domain of mode 1, do not prevent these intersections (e.g., it could lead to chattering \cite[Section~1.2.4]{Liberzon2012} around the switching surface, because any trajectory starting in mode 2 will cross the switching surface to mode 1, by the stability of $x^*$ in $\SD$).

The following assumption handles these scenarios.
\cite[Figure 6.7]{Gomes2019d} illustrates the concepts in this assumption and shows the impact of different \PCRay{s}.

\begin{assumption}
  \label{ass:perppc}
  Given a planar \hsys{} $\SD_H$, let
  \begin{align}
    \notag
    \GoodSwitchingSurface & = \{\, x \in \SwitchingSurface \mid \nabla \swfun(x) \cdot f_2(x) < 0 \,\}, \\
    \label{eq:T2}
    \TimeModeTwoSp & = \inf_{x \in \GoodSwitchingSurface} \rtime{f_2}{\PointCarreLine}{x}, \text{ and}\\
    \notag
    \GoodD & = \{\, \Flowi{2}{x}{h} \mid x \in \GoodSwitchingSurface, 0 \le h < \TimeModeTwoSp \,\}.
  \end{align}
  We assume that
  for every point $x \in \GoodD$,
  $\rtime{f_1}{\PointCarreLine}{x}$ is finite and smaller than $\rtime{f_1}{\SwitchingSurface}{x}$ (which may be infinite).
\end{assumption}

Under \cref{ass:perppc}, and a sufficiently small $H$ in $\SD_H$, trajectories that start in $\GoodPointCarreLine$ acquire a predictable pattern, illustrated in \cite[Figure 6.8]{Gomes2019d}.
We now define the notation and restrictions to represent these trajectories.

\begin{figure}[htb]
	\centering
	\begin{subfigure}[b]{0.25\textwidth}
	\centering
	\includegraphics[width=0.9\textwidth]{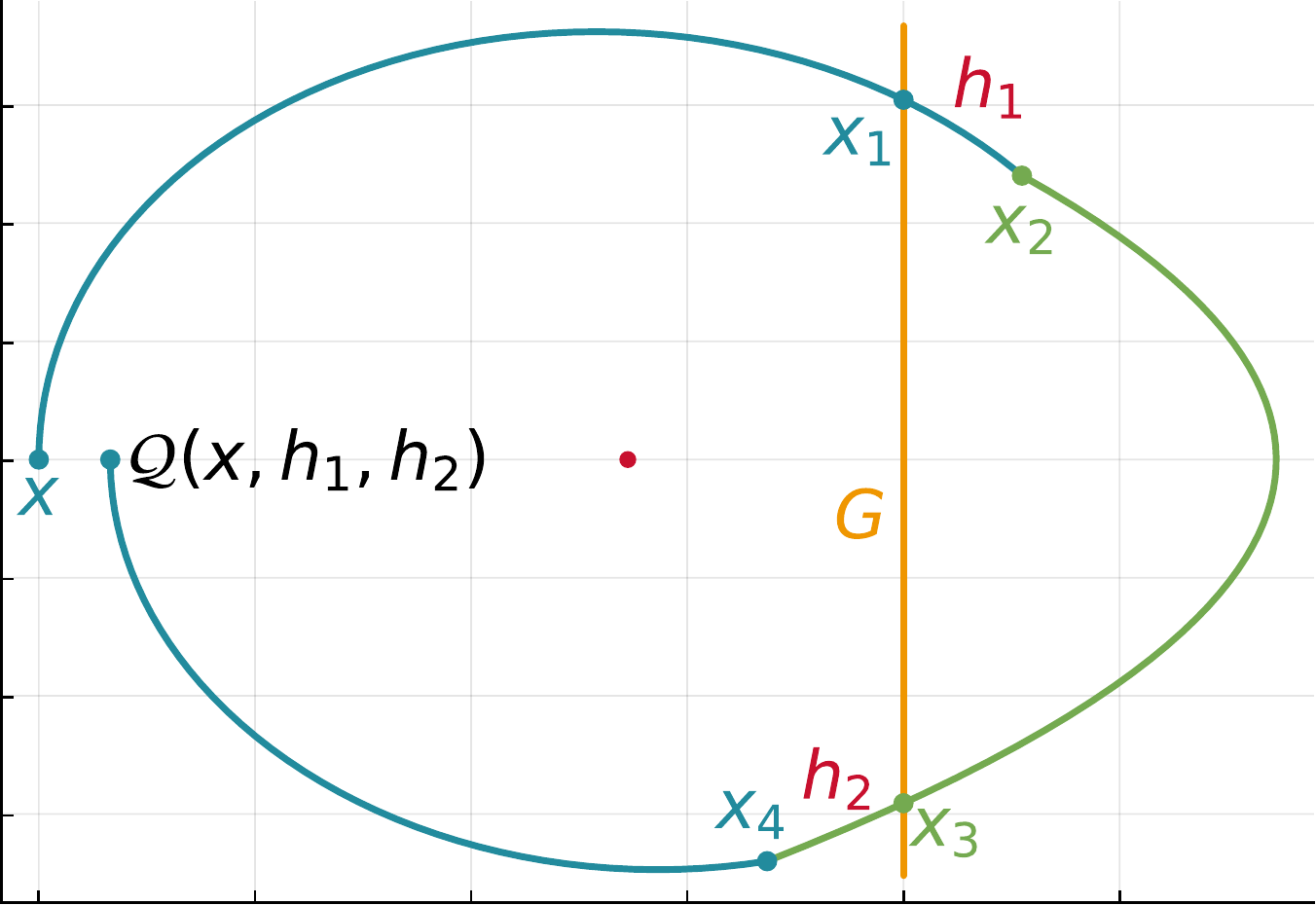}
	\caption{
	Illustration of $\hpcmap{h_1, h_2}{x}$.
	}
	\label{fig:hpcmap}
	\end{subfigure}%
	~%
	\begin{subfigure}[b]{0.25\textwidth}
	\centering
	\includegraphics[width=0.9\textwidth]{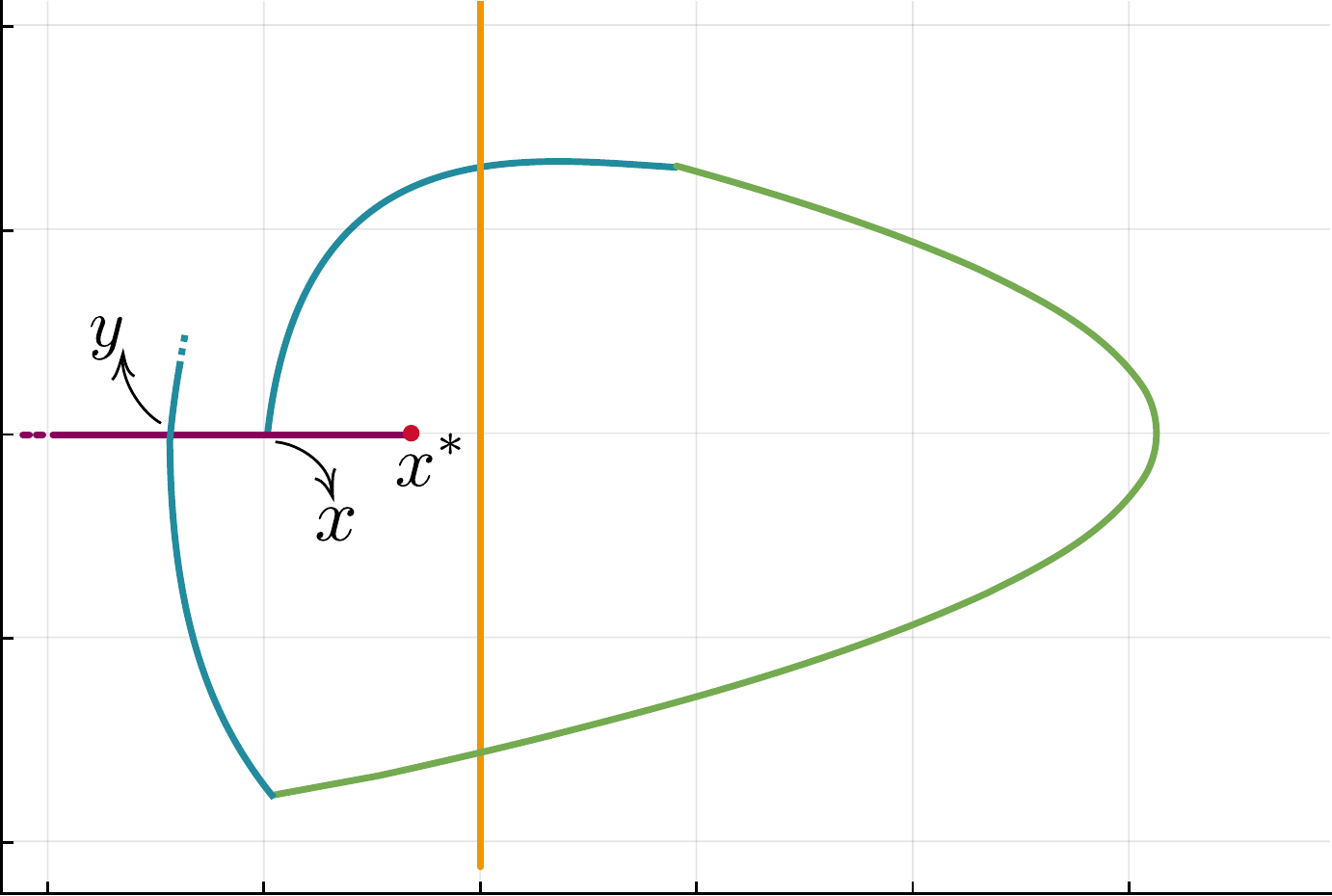}
	\caption{Illustration of \cref{lem:U}.}
	\label{fig:example_lemmaU}
	\end{subfigure}
	\caption{Illustration of the mapping $\hpcmap{h_1, h_2}{x}$, and \cref{lem:U}, for the bouncing ball example (\cref{sec:motivating_example}). The horizontal axis refers to position, and vertical to velocity.}
	\end{figure}

Given a point $x_0 \in \GoodPointCarreLine$ and delays $h_1, h_2 \ge 0$, let
\begin{align}
  \label{eq:x1x2}
  x_1 & = \pcmap{f_1}{\SwitchingSurface}{x_0} &
  x_2 & = \Flowf{f_1}{x_1}{h_1} \\
  \label{eq:x3x4}
  x_3 & = \pcmap{f_2}{\SwitchingSurface}{x_2} &
  x_4 & = \Flowf{f_2}{x_3}{h_2}.
\end{align}
If $h_1 < \rtime{f_1}{\SwitchingSurface}{x_1}$ and $h_2 < \TimeModeTwoSp$,
by \cref{ass:perppc},
the following \exec{} is admitted by planar $\SD_{\max(h_1, h_2)}$.
\begin{align}\label{eq:tq_over_pccurve}
  \hmapexec{h_1, h_2}{x_0} &
  = 
        \joinexec{
          \joinexec{
            \execpcmap{1}{\SwitchingSurface}{x_0}
          }{
            \execFlow{1}{h_1}{x_1}
          }
        }{
          \execpcmap{2}{\SwitchingSurface}{x_2}
        }\\ &
  \quad \joinexec{
      \joinexec{
      }{
        \execFlow{2}{h_2}{x_3}
      }
    }{
      \execpcmap{1}{\PointCarreLine}{x_4}
    }.
\end{align}
We also introduce the notation
\begin{equation*}
\hpcmap{h_1, h_2}{x_0} \triangleq \pcmap{f_1}{\PointCarreLine}{x_4},
\end{equation*}
illustrated in \cref{fig:hpcmap}. 

The following lemma is illustrated in \cref{fig:example_lemmaU}. 

\begin{lemma}
  \label{lem:U}
  Given a delay $H < \TimeModeTwoSp$ (where $\TimeModeTwoSp$ is defined in \ceqref{eq:T2}), a planar \hsys{} $\SD_H$ as given in \cref{def:SDH}, with a \PCRay{} $\PointCarreLine$ as defined in \ceqref{eq:pc}.
  Under \cref{ass:closedinj,ass:perppc},
  if the equilibrium $x^*$ is not GAS for $\SD_H$, and $\SD_H$ does not admit any closed orbit,
  then there exists two points $x, y \in \GoodPointCarreLine$ and delays $0 \le h_1, h_2 \le H$ such that
  $y = \hpcmap{h_1, h_2}{x}$ and $\enorm{y} > \enorm{x}$.
\end{lemma}

The following lemma shows that, by adjusting the delay in $\SD_H$, one can find admissible trajectories whose successive intersections with the \PCRay{} are getting closer to the equilibrium.
The lemma statement is illustrated in \cite[Figure 6.11]{Gomes2019d}.

\begin{lemma}
  \label{lem:SH}
  Given a delay $H < \TimeModeTwoSp$ (where $\TimeModeTwoSp$ is defined in \ceqref{eq:T2}), a planar \hsys{} $\SD_H$ as given in \cref{def:SDH}, with a \PCRay{} $\PointCarreLine$ as defined in \ceqref{eq:pc}.
  Under \cref{ass:closedinj},
  if there exist two points $x, y \in \GoodPointCarreLine$ and delays $0 \le h_1, h_2 \le H$ such that
  $y = \hpcmap{h_1, h_2}{x}$ and $\enorm{y} > \enorm{x}$
  then there exist
  delays $0 \le \overline{h_1}, \underline{h_1} \le h_1$ and $0 \le \overline{h_2}, \underline{h_2} \le h_2$ such that
  \[ x \in [\hpcmap{\underline{h_1}, \underline{h_2}}{x}, \hpcmap{\overline{h_1}, \overline{h_2}}{x}]_{\PointCarreLine}. \]
\end{lemma}

From the existence of an unstable trajectory,
\cref{lem:U} combined with \cref{lem:SH}
ensures the existence of a point $x$ such that for different delays,
the next intersection with the \PCRay{} is either ``above'' or ``below'' it.
We now prove the continuity of the \PCMap{} to show that there are delays such
that the next intersection is exactly $x$.


The continuity of \PCMap{} is classically deduced from the \emph{Implicit Function Theorem}, see e.g. \cite[Section~5.2, Theorem~2.1]{Robinson1998}.

\begin{lemma}[Continuity of \PCMap{}]
  \label{lem:cont}
  Consider two continuously differentiable functions $f, g : \setreal^2 \to \setreal$, and
  let $G = \{\, x \mid g(x) = 0\,\}$.
  For any open set $A$ such that $A \cap G = \emptyset$:
  if for all $x \in A$, $z = \pcmap{f}{G}{x}$ is defined,
  and $\nabla g(z) \cdot f(z) \neq 0$,
  then the \PCMap{} $\pcmap{f}{G}{x}$ restricted to $A$ is continuous.
\end{lemma}

\begin{remark}\label{rem:continuity}
  The continuity of the \PCMap{} cannot be readily generalized to the
  \hsys{} context.
  Given planar $\SD_H$,
  the \PCMap{} may have discontinuities in $x$ and in the \swdl{}.
  Discontinuities in $x$ may happen if the \exec{} is tangent to
  the switching surface
  while discontinuities in the \swdl{} may happen when the time spent in
  $\Dom(2)$ is exactly the \swdl{}.
  \cite[Figure 6.13]{Gomes2019d} illustrates these examples.

  While the \PCMap{} is not continuous everywhere, under some assumptions
  that will be stated, our argument only relies on its continuity on regions
  where it is locally continuous.
\end{remark}

In view of the statement of \cref{lem:cont},
we add the following two assumptions.

\begin{assumption}
  \label{ass:perpsw}
  We assume that for every point $x \in \interior(\Dom(2))$,
  the point $y = \pcmap{f_2}{\SwitchingSurface}{x}$ is such that
  the surface normal $\nabla \swfun(y)$ and $f_2(y)$ are not perpendicular,
  that is, $\nabla \swfun(y) \cdot f_2(y) \neq 0$.
\end{assumption}


\begin{assumption}
  \label{ass:perpsw1}
  We assume that for every point $x \in \GoodPointCarreLine$,
  the point $y = \pcmap{f_1}{\SwitchingSurface}{x}$ is such that
  the surface normal $\nabla \swfun(y)$ and $f_1(y)$ are not perpendicular,
  that is, $\nabla \swfun(y) \cdot f_1(y) \neq 0$.
\end{assumption}

\begin{lemma}
  \label{poincare_map_continuous}
  Consider a planar BMHA $\SD$ as given in \cref{def:bms} satisfying \cref{ass:perpsw},
  and a \PCRay{} $\PointCarreLine$ satisfying \cref{ass:perppc,ass:perpsw1}.
  For any point $x_0 \in \GoodPointCarreLine$,
  the \PCMap{} $\hpcmap{h_1, h_2}{x_0}$ is defined and continuous with respect
  to $x_0$, $h_1 \in [0; \rtime{f_1}{\SwitchingSurface}{\pcmap{f_1}{\SwitchingSurface}{x}}[$ and $h_2 \in [0; \TimeModeTwoSp[$ where $\TimeModeTwoSp$ is defined in \ceqref{eq:T2}.
\end{lemma}

\begin{proposition}\label{claim:instability_implies_orbit}
  Given a delay $H < \TimeModeTwoSp$ (where $\TimeModeTwoSp$ is defined in \ceqref{eq:T2}), a planar \hsys{} $\SD_H$ as given in \cref{def:SDH}, with a \PCRay{} $\PointCarreLine$ as defined in \ceqref{eq:pc}.
  Under \cref{ass:closedinj,ass:perppc,ass:perpsw,ass:perpsw1},
  if the equilibrium $x^* \in \setreal^2$ of $\SD_H$ is not GAS and $H < \TimeModeTwoSp$,
  then $\SD_H$ admits a closed orbit.
\end{proposition}

The above claim allows us to solve \cref{eq:general_problem} by solving the following problem:
\begin{aligneq}\label{eq:better_problem}
\hat{\sigma}(\SD) = \inf_{H \geq 0} H \text{ s.t. } \SD_H \text{ admits a closed orbit}.
\end{aligneq}

\begin{theorem}\label{orbits_equivalent_stability}
  Given a planar BMHA $\SD$ as given in \cref{def:bms} satisfying \cref{ass:perpsw}
  and a \PCRay{} $\PointCarreLine$ as defined in \ceqref{eq:pc}
  satisfying \cref{ass:closedinj,ass:perppc,ass:perpsw1},
  and $\sigma(\SD) < \TimeModeTwoSp$ (defined in \ceqref{eq:T2}),
  the identity $\sigma(\SD) = \hat{\sigma}(\SD)$ holds. 
\end{theorem}
\ifproofs %
\begin{proof}
    The inequality $\sigma(\SD) \leq \hat{\sigma}(\SD)$ follows from \cref{claim:orbit}
    and we show $\sigma(\SD) \geq \hat{\sigma}(\SD)$ by contradiction.
    If $\sigma(\SD) < \hat{\sigma}(\SD)$ then there exists a delay $H$ such that $\sigma(\SD) < H < \hat{\sigma}(\SD)$.
    Since $\sigma(\SD) < H$, the equilibrium $x^*$ is not GAS for $\SD_H$.
    Therefore, by \cref{claim:instability_implies_orbit}, $\SD_H$ admits a closed orbit.
    This is in contradiction with $H < \hat{\sigma}(\SD)$.
\end{proof}
\else %
\fi%

%% file: application.tex
\section{Results}
\label{sec:result}

This section shows how \cref{orbits_equivalent_stability} can be applied to compute the \MSD{} $\sigma(\SD)$ in \cref{eq:general_problem}.
Then, it illustrates the solution procedure to the bouncing ball example, introduced in \cref{sec:motivating_example}.

As discussed in \cref{sec:orbit},
we can restrict our attention to orbits of the form $\hmapexec{h_1, h_2}{x}$ where
$x \in \GoodPointCarreLine$.
That is, we have
\begin{multline}
  \hat{\sigma}(\SD) = \inf_{x \in \GoodPointCarreLine, h_1, h_2}
  \{\,
    \max(h_1, h_2) \mid
    x \in \GoodPointCarreLine,\\
    h_1, h_2 \ge 0,
    \hpcmap{h_1, h_2}{x} = x
  \,\}.
\end{multline}
This constrained 3-dimensional nonlinear optimization problem can be reduced to
the following unconstrained 2-dimensional nonlinear optimization problem:
\begin{equation}
  \label{eq:2Dopt}
  \hat{\sigma}(\SD) = \inf_{(x, h_1)} \max(h_1, \opt{h_2}(x, h_1))
\end{equation}
where $\opt{h_2}(x, h_1) = \min\{\, h_2 \mid \hpcmap{h_1, h_2}{x} = x \,\}$.
This reduction is possible because,
given fixed values of $x$ and $h_1$, the value of $\opt{h_2}(x, h_1)$ is straightforward to compute.
Indeed, let $y \in \SwitchingSurface$ be such that $\pcmap{f_1}{\PointCarreLine}{y} = x$,
we have
\[ \opt{h_2}(x, h_1) = \rtime{f_2}{\execvar}{x_2} - \rtime{f_2}{\SwitchingSurface}{x_2} \]
where $\execvar = \execpcmap{1}{\PointCarreLine}{y}$ and $x_2$ is defined in \ceqref{eq:x1x2}.
For a given accuracy,
$\opt{h_2}(x, h_1)$ can be computed using classical simulation methods for nonlinear systems (see e.g. \cite{Cellier2006}) as follows.
For a given point $x \in \GoodPointCarreLine$,
we precompute the trajectory $\execvar$ with a time steps determined by the required accuracy.
Then we compute the point $\rtime{f_2}{\execvar}{x_2}$ by first simulating
the trajectory with a coarse time step.
Let $t, z$ be the last element of the sampled trajectory before the intersection
with $\execvar$.
We have $\rtime{f_2}{\execvar}{x_2} = t + \rtime{f_2}{\execvar}{z}$
and $\rtime{f_2}{\execvar}{z}$ is smaller than the time step used to
simulate the trajectory starting at $x_2$.
We can therefore estimate $\rtime{f_2}{\execvar}{z}$ with a refined time step.
This procedure can be applied recursively.

In the example introduced in \cref{sec:motivating_example},
the optimal solution of the problem in \ceqref{eq:2Dopt} with accuracy $10^{-9}$
is given at $\opt{x} \approx 0.24579453$, $\opt{h_1} = 0$ and $h_2 = \opt{h_2}(\opt{x}, 0) \approx 0.0014128697$.
We illustrate the objective function along the line $h_1 = 0$ in \fighplot{}.
We used \cite{rackauckas2017differentialequations} to simulate the nonlinear
system.

\cite[Section 6.5.1]{Gomes2019d} provides a comparison of our algorithm with a trial-and-error approach using SpaceEx.

%% file: conclusion.tex
\section{Conclusion}
\label{sec:conclusion}

Motivated by practical problems in the simulation of hybrid systems, we have studied how a delay in the detection of mode transitions can impact the quality of the result of the numerical simulation.  It turns out that this delay may have a crucial impact on the result, as it may turn a stable hybrid system into an unstable behaviour in the numerics.

Our goal was to study this phenomenon, and we have focused here on planar systems.
A natural first research question aiming at understanding the problem is to characterize the threshold between stability and instability of the simulated trajectories, when the true system is stable.  Already for this simple case, it can be hard to compute the maximal allowed delay ensuring stability of the trajectories.

We have used classical techniques in the analysis of dynamical systems, such as Poincaré maps and topological arguments.  However, we showed that in hybrid systems, more complex phenomena can occur, which make these classical tools insufficient to solve the problem.  We pushed further these techniques, which allowed us, under mild assumptions, to compute this maximal delay for planar systems.
Moreover, we showed that traditional reachability analysis techniques do not scale to solve this problem satisfactorily.

We hope that this work will provide a better understanding of the problem, of high importance in numerical simulation, and will lead to the estimation of the maximal allowed delay for more complex, or higher dimensional, hybrid systems than the ones studied here.

%% file: proofs.tex
\section{Proofs}

\begin{proof}[\cref{lem:intersect}]
  Let $\setvar$ be the path component of $\mathcal{X} \setminus (\statetimeexec{t_1}{t_2} \cup [x_2, x_1]_{\PointCarreLine})$, that contains $x^* \in \setvar$.
  By \Cref{ass:closedinj},
  there is a path in $\closure{\PointCarreLine}$ from $x^*$ to $x_2$.
  This path is contained in $\closure{\setvar}$ as it cannot intersect the interior
  of $\statetimeexec{t_1}{t_2}$ by definition of $t_2$.
  Therefore, by the injectiviy assumption of \Cref{ass:closedinj},
  the set $\setvar$ must contain all $x \in \PointCarreLine$ such
  that $\enorm{x} \le \enorm{x_2}$.

  As $\nabla \pcfun(x) \cdot f_1(x) > 0$ for all $x \in \PointCarreLine$ (see \ceqref{eq:pc})
  and $\enorm{x_2} \ge \enorm{x_1}$, 
  there exists $\epsilon > 0$ such that $\statetimeexec{t_2}{t_2 + \epsilon} \cap U = \emptyset$.
  As the \exec{} at $t_2 + \epsilon$ is not in $\setvar$,
  if $\statetimeexec{t_2}{\infty} \cap \setvar \neq \emptyset$ then
  $(q_3, t_3, x_3) \in \fullpcmap{\execvar}{\closure{U}}{t_2}$ is defined.
  The point $x_3$ cannot belong to $[y, x]_{\PointCarreLine}$ as
  $\nabla \pcfun(x) \cdot f_1(x) > 0$ hence
  $x_3 \in \statetimeexec{t_1}{t_2}$.
  For any $(q, t, x_3) \in \timeexec{t_1}{t_2}$, we must have $q \neq q_3$ since
  by \cite[Lemma 1]{Gomes2019d}, two \execs{} in the same mode cannot intersect.
\end{proof}

\begin{proof}[\cref{lem:invariant_mode_1}]
  Let $\setvar$ be the path component of $\mathcal{X} \setminus (\statetimeexec{t_1}{t_2} \cup [x_2, x_1]_{\PointCarreLine})$
  containing the equilibrium $x^*$.
  If $\statetimeexec{t_2}{\infty} \not\subseteq \closure{U}$ then
  given
  $$(q', t', x') \in \fullpcmap{\execvar}{\statetimeexec{t_1}{t_2} \cup [x_2, x_1]_{\PointCarreLine}}{t_2},$$
  we know that $x' \in \statetimeexec{t_1}{t_2}$ since $\nabla \pcfun(x) \cdot f_1(x) > 0$ for all $x \in \PointCarreLine$.
  As $\statetimeexec{t_1}{t_2} \subseteq \Dom(1)$, $q' = 1$ which
  is impossible by \cite[Lemma 1]{Gomes2019d}.
  Therefore $\statetimeexec{t_2}{\infty} \subseteq \closure{U} \subseteq \Dom(1)$.
\end{proof}

\begin{proof}[\cref{lem:U}]
  Since $\SD_H$ is not GAS, it admits \anexec{} $\execvar$ that does not converge to
  the equilibrium $x^*$.

  The proof is divided into two parts: first we prove that \exec{} $\execvar$ contains infinitely many mode transitions.
  Then we use that fact to prove that successive intersections satisfy the claim in the lemma.
  In the second part, we make use of a result that is proved later (\cref{poincare_map_continuous}). \Cref{poincare_map_continuous} does not depend on \cref{lem:U}.

  We now prove by contradiction that \exec{} $\execvar$ contains infinitely many mode transitions.
  If it contains finitely many transitions, there is a mode $q_\infty$ and time $t_\infty$ such that,
  for all $(q, t, x) \in \execvar$ with $t > t_\infty$, we have $q = q_\infty$.
  We cannot have $q_\infty = 1$ as the equilibrium $x^*$ is GAS for mode 1.
  Similarly, we cannot have $q_\infty = 2$ since the equilibrium $x^*$ is GAS for $\SD$.
  In this case,
  $(2, t', x') \in \fullpcmap{\execvar}{\SwitchingSurface}{t_\infty}$ is defined.
  Let $(2, t'', x'') \in \fullpcmap{\execvar}{\SwitchingSurface}{t'}$.
  Then, by \Cref{ass:perppc}, there must be $t''', x'''$ such that $(1, t''', x''') \in \timeexec{t'}{t''}$.
  This concludes the proof that \exec{} $\execvar$ contains infinitely many mode transitions.

  Now we focus on the second part of the proof.
  By \Cref{ass:perppc}, after each transition from mode 2 to mode 1,
  the \exec{} $\execvar$ must intersect the \PCRay{} before it
  transitions from mode 1 to mode 2.
  Therefore the \exec{} intersects the \PCRay{} an infinite amount of times
  and is in mode 1 each time it intersects it.
  Let $(1, t_1, x_1), (1, t_2, x_2), \ldots \in \execvar$ with $t_1 < t_2 < \cdots$
  and $x_1, x_2, \ldots \in \PointCarreLine$ be this sequence intersections.

  There are no $i \neq j$ such that $\enorm{x_i} = \enorm{x_j}$
  because, by \Cref{ass:closedinj}, $\enorm{x_i} = \enorm{x_j}$ implies that
  $x_i = x_j$ which contradicts the statement that no closed orbit is admitted.
  There are two possible cases.
  \begin{compactitem}
    \item
    The sequence $(\enorm{x_i})_i$ is decreasing with $i$.
    Since it is bounded below by 0, it converges.
    Let $n = \lim_{i \to \infty} \enorm{x_i}$.
    Since the \exec{} is not stable, $n > 0$.
    By \Cref{ass:closedinj}, there exists a unique $x \in \PointCarreLine$
    such that $\enorm{x} = n$.
    By \cref{poincare_map_continuous}, $\SD_H$ admits a closed orbit containing $x$
    which contradicts the statement. Hence this case cannot be true.
    \item
    There exists a $k$ such that $\enorm{x_{k+1}} > \enorm{x_k}$.
    Let
    $$K = \{\, k \mid \enorm{x_{k+1}} > \enorm{x_k} \,\}.$$
    There are two possible sub-cases.
    \begin{compactitem}
      \item If there exists $k \in K, t \in \setreal, x \in \mathcal{X}$ such that
      $(2, t, x) \in \timeexec{t_k}{t_{k+1}}$ then we are done.
      \item Otherwise, for all $k \in K$, $x_{k+1} = \pcmap{f_1}{\PointCarreLine}{x_k}$.
      Note that, by \cref{coro:decreasing}, it must be the case that the trajectory between $x_k$ and $x_{k+1}$ is going into the domain of mode 2.
      Let $j \in K$
      be such that $j+1 \notin K$, and
      $$(\tilde{q}, \tilde{t}, \tilde{x}) \in \fullpcmap{\execvar}{\statetimeexec{t_j}{t_{j+1}}}{t_{j+1}},$$
      as illustrated in \cite[Figure 6.10]{Gomes2019d}.
      By \cref{lem:intersect}, $\tilde{q} = 2$.
      Let $\bar{t}$ be such that $(1, \bar{t}, \tilde{x}) \in \timeexec{t_j}{t_{j+1}}$.
      If $\tilde{x} \in \Dom(1)$, as in \cite[Figure 6.10]{Gomes2019d}, then $\timeexec{\bar{t}}{\tilde{t}}$ is
      a closed orbit which contradicts the statement.
      Otherwise, $\tilde{x} \in \Dom(2)$\footnote{In this case, $\timeexec{\bar{t}}{\tilde{t}}$ is not a closed orbit since it requires a transition between mode 2 and mode 1 at $\tilde{t}$ which is not possible since $\bar{x} \in \Dom(2)$.}, as shown in \cite[Figure 6.10]{Gomes2019d}.
      Let $\setvar$ be the path component of $\mathcal{X} \setminus \statetimeexec{\bar{t}}{\tilde{t}}$ 
      containing the points $x \in \PointCarreLine$ such that $\enorm{x} < \enorm{x_{j+1}}$
      and
      $j' = \min \{\, k \in K \mid k > j+1, \enorm{x_{k+1}} > \enorm{x_{j+1}} \,\}$.
      Since $x_{j'+1} \notin \closure{\setvar}$, by \cite[Lemma 1]{Gomes2019d},
      $(q', t', x') \in \fullpcmap{\execvar}{\statetimeexec{\bar{t}}{\tilde{t}}}{t_{j'}}$ is defined and is such that
      $t_{j'} \le \bar{t} \le t_{j'+1}$,
      $q' = 1$ and there exists
      $(2, t'', x') \in \timeexec{\bar{t}}{\tilde{t}}$.
      The \exec{} $\timeexec{t''}{t'}$ is a closed orbit which contradicts
      the statement.
    \end{compactitem}
  \end{compactitem}
\end{proof}

\begin{figure}[htb]
  \begin{center}
    \includegraphics[width=\defaultImageWidth\textwidth]{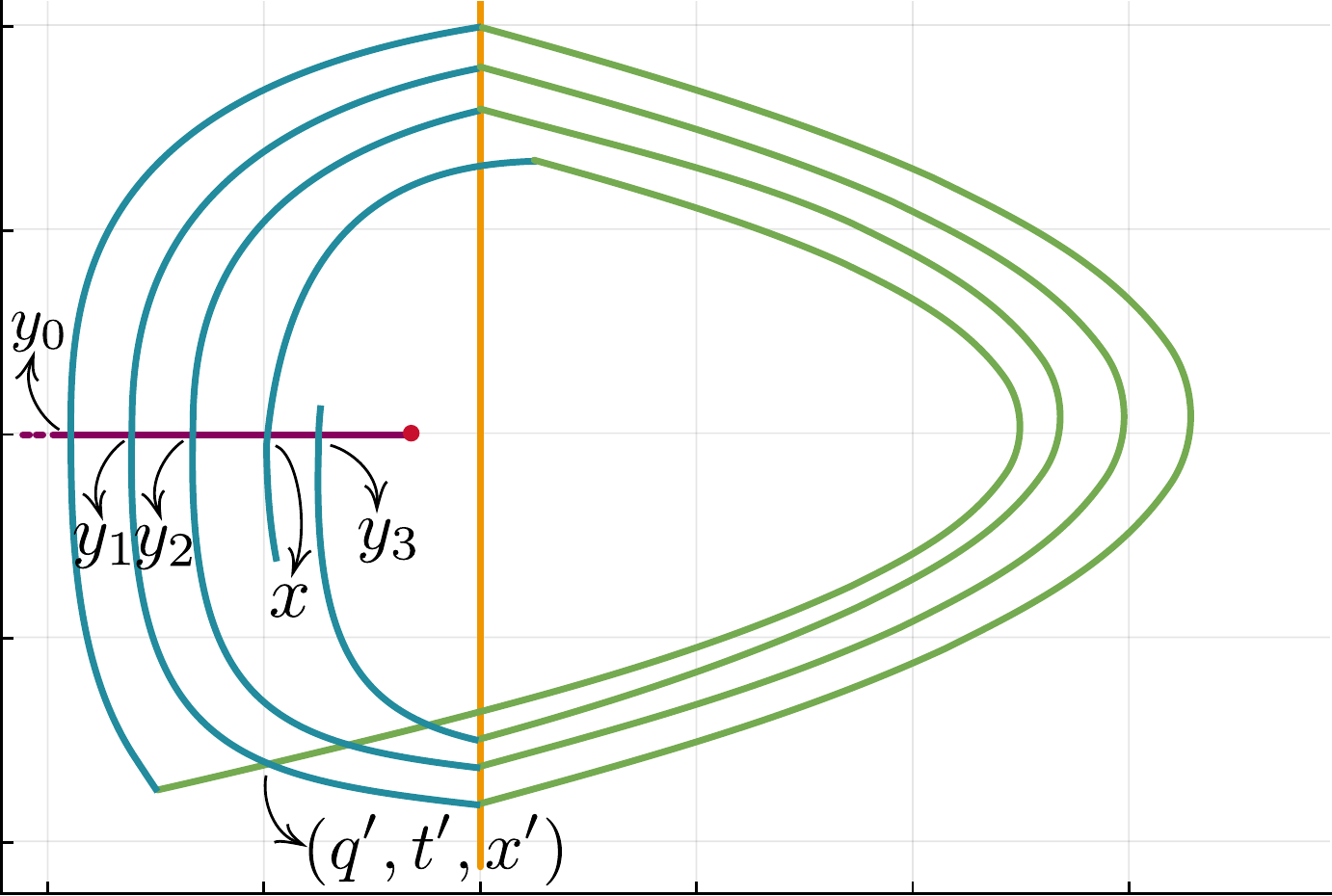}
    \caption{Illustration of the sequence defined in \ceqref{eq:sequence}.}
    \label{fig:lema_sh}
  \end{center}
\end{figure}

\begin{proof}[\cref{lem:SH}]
  Consider \anexec{} $\execvar$ of $\SD$ starting at $y$.
  Let $y_0 = y$, $t_0 = 0$ and
  \begin{equation}\label{eq:sequence}
  (1, t_{k+1}, y_{k+1}) \in \fullpcmap{\execvar}{\PointCarreLine}{t_k},
  \end{equation}
  for $k = 0, 1, \ldots$ 
  This sequence is illustrated in \cref{fig:lema_sh}.
  By \cref{coro:decreasing}, the sequence $(\enorm{y_k})_k$ is non-increasing in $k$.
  Let $j = \inf \{\, k \mid \enorm{y_k} < \enorm{x} \,\}$,
  we know that $j < \infty$ since the equilibrium $x^*$ is GAS for $\SD$.
  In \cref{fig:lema_sh}, $j=3$.
  We show by recurrence the existence of
  a sequence of delays $h_{1, k} \le h_1$ and $h_{2, k} \le h_2$ for $k = 0, 1, \ldots, j$
  such that $y_{k} = \hpcmap{h_{1, k}, h_{2, k}}{x}$, where $y_{k}$ is defined in \ceqref{eq:sequence}.
  We set $h_{1, 0} = h_1$ and $h_{2, 0} = h_2$.
  For $k \ge 0$, given $h_{1, k} \le h_1$ and $h_{1, k} \le h_2$,
  let
  $$ (q', t', x') \in \fullpcmap{\execvar}{\statehmapexec{h_{1, k}, h_{2, k}}{x}}{t_k}, $$
  as illustrated in \cref{fig:lema_sh},
  and
  $$(q'', t'', x') \in \hmapexec{h_{1, k}, h_{2, k}}{x},$$
  where $\hmapexec{h_{1, k}, h_{2, k}}{x}$ is defined in \ceqref{eq:tq_over_pccurve}.
  By \cref{lem:intersect}, we only have the following two cases to consider%
  \footnote{Note that in both cases, we have
    $\hmapexec{h_{1, k+1}, h_{2, k+1}}{x} =
    \joinexec
    {\timeslice{\hmapexec{h_{1, k}, h_{2, k}}{x}}{0}{t''}}
    {\timeexec{t'}{t_{k+1}}}$.}:
  if $q' = 2$, $q'' = 1$ and $x' \in \Dom(2)$ then we set
  $h_{2, k+1} = 0$ and $h_{1, k+1} = t'' - \rtime{f_1}{\SwitchingSurface}{x}$;
  otherwise, if $q' = 1$, $q'' = 2$ and $x' \in \Dom(1)$ then we set $h_{1, k+1} = h_{1, k}$
  and $h_{2, k+1} = t'' - (\rtime{f_1}{\SwitchingSurface}{x} + h_{1, k} + \rtime{f_2}{\SwitchingSurface}{\Flowf{f_1}{\pcmap{f_1}{\SwitchingSurface}{x}}{h_{1,k}}})$.
  The lemma is proved with
  $\overline{h_1} = h_{1, j-1}$, $\underline{h_1} = h_{1, j}$,
  $\overline{h_2} = h_{2, j-1}$and $\underline{h_2} = h_{2, j}$.
\end{proof}

\begin{proof}[\cref{lem:cont}]
  Consider a point $x_0 \in A$.
  Let $t_0 = \rtime{f}{G}{x_0}$
  and $F(x, t) = g(\Flowf{f}{x}{t})$,
  we have $F(x_0, t_0) = 0$.
  Furthermore
  \[ \frac{\partial F}{\partial t}(x_0, t_0) = \nabla g(\Flowf{f}{x}{t}) \cdot \frac{\partial \Flowf{f}{x}{t}}{\partial t} = \nabla g(\Flowf{f}{x}{t}) \cdot f(\Flowf{f}{x}{t}) \]
  which is nonzero by assumption.
  Therefore, by the Implicit Function Theorem \cite[Section~5.2, Theorem~2.1]{Robinson1998}, there exists an open set $U \subseteq \setreal^2$
  containing $x_0$
  and a unique continuously differentiable function $h$ such that $F(x, h(x)) = 0$ for each $x \in U$.
  By the uniqueness of $h$, we know that $\pcmap{f}{G}{x} = \Flowf{f}{x}{h(x)}$ for each $x \in U$.
  Since the function $h(x)$ is continuous in $x$ and the function $\Floww$ is continuous in $x$ and $t$ by \cite[Lemma 1]{Gomes2019d},
  the \PCMap{} $\pcmap{f}{G}{x}$ is continuous in $x_0$.
\end{proof}

\begin{proof}[\cref{poincare_map_continuous}]
  Let $x_1, x_2, x_3, x_4$ be as defined in \ceqreff{eq:x1x2}{eq:x3x4}
  and $x_5$ denote $\hpcmap{h_1, h_2}{x}$.

  By \Cref{ass:perppc} and \cref{lem:cont}, $x_5$ depends continuously on $x_4$ and
  by \cite[Lemma 1]{Gomes2019d}, $x_4$ depends continuously on $h_2$
  hence $x_5$ depends continuously on $h_2$.

  By \cite[Lemma 1]{Gomes2019d}, $x_4$ depends continuously on $x_3$
  and by \Cref{ass:perpsw} and \cref{lem:cont}, $x_3$ depends continuously on $x_2$
  hence $x_5$ depends continuously on $x_2$.

  By \cite[Lemma 1]{Gomes2019d}, $x_2$ depends continuously on $h_1$ hence
  $x_5$ depends continuously on $h_1$.

  By \cite[Lemma 1]{Gomes2019d}, $x_2$ depends continuously on $x_1$ and
  by \Cref{ass:perpsw1} and \cref{lem:cont}, $x_1$ depends continuously on $x_0$
  hence $x_5$ depends continuously on $x_0$.
\end{proof}

\begin{proof}[\cref{claim:instability_implies_orbit}]
  Suppose by contradiction that there exists $H$ such that $\SD_H$ is not GAS but does not admit any closed orbit.
  Consider one such value of $H$.
  By \cref{lem:U} and \cref{lem:SH}, there a point $x \in \GoodPointCarreLine$,
  delays $0 \le \overline{h_1}, \underline{h_1} \le H$ and $0 \le \overline{h_2}, \underline{h_2} \le H$
  such that
  \[ x \in [\hpcmap{\underline{h_1}, \underline{h_2}}{x}, \hpcmap{\overline{h_1}, \overline{h_2}}{x}]_{\PointCarreLine}. \]
  By \cref{poincare_map_continuous}, there exists
  $\min(\overline{h_1}, \underline{h_1}) \le h_1 \le \max(\overline{h_1}, \underline{h_1})$ and
  $\min(\overline{h_2}, \underline{h_2}) \le h_2 \le \max(\overline{h_2}, \underline{h_2})$ such that
  $\hpcmap{h_1, h_2}{x} = x$.
  This contradicts the fact that $\SD_H$ admits no closed orbit. Hence the proof is complete.
\end{proof}

%% file: root.bbl
\ifdefined\DeclarePrefChars\DeclarePrefChars{'’-}\else\fi
\begin{thebibliography}{10}
\providecommand{\url}[1]{#1}
\csname url@samestyle\endcsname
\providecommand{\newblock}{\relax}
\providecommand{\bibinfo}[2]{#2}
\providecommand{\BIBentrySTDinterwordspacing}{\spaceskip=0pt\relax}
\providecommand{\BIBentryALTinterwordstretchfactor}{4}
\providecommand{\BIBentryALTinterwordspacing}{\spaceskip=\fontdimen2\font plus
\BIBentryALTinterwordstretchfactor\fontdimen3\font minus
  \fontdimen4\font\relax}
\providecommand{\BIBforeignlanguage}[2]{{%
\expandafter\ifx\csname l@#1\endcsname\relax
\typeout{** WARNING: IEEEtran.bst: No hyphenation pattern has been}%
\typeout{** loaded for the language `#1'. Using the pattern for}%
\typeout{** the default language instead.}%
\else
\language=\csname l@#1\endcsname
\fi
#2}}
\providecommand{\BIBdecl}{\relax}
\BIBdecl

\bibitem{Jungers2009}
R.~Jungers, \emph{The Joint Spectral Radius: Theory and Applications}.\hskip
  1em plus 0.5em minus 0.4em\relax {Springer Science \& Business Media}, 2009,
  vol. 385.

\bibitem{Bauer2012}
N.~W. Bauer, P.~J.~H. Maas, and W.~P. M.~H. Heemels, ``Stability analysis of
  networked control systems: {{A}} sum of squares approach,''
  \emph{Automatica}, vol.~48, no.~8, pp. 1514--1524, 2012.

\bibitem{Conti2016}
\BIBentryALTinterwordspacing
F.~Conti and O.~Khatib, ``A {{Framework}} for {{Real}}-{{Time Multi}}-{{Contact
  Multi}}-{{Body Dynamic Simulation}},'' in \emph{Robotics {{Research}}}.\hskip
  1em plus 0.5em minus 0.4em\relax Cham: {Springer International Publishing},
  2016, vol. 114, pp. 271--287. [Online]. Available:
  \url{http://link.springer.com/10.1007/978-3-319-28872-7\_16}
\BIBentrySTDinterwordspacing

\bibitem{Gomes2019d}
\BIBentryALTinterwordspacing
C.~Gomes, ``Property preservation in co-simulation.'' [Online]. Available:
  \url{https://repository.uantwerpen.be/docman/irua/57f437/163840.pdf}
\BIBentrySTDinterwordspacing

\bibitem{Jungers2012}
R.~M. Jungers, A.~D'Innocenzo, and M.~D. Di~Benedetto, ``Feedback stabilization
  of dynamical systems with switched delays,'' in \emph{2012 {{IEEE}} 51st
  {{IEEE Conference}} on {{Decision}} and {{Control}} ({{CDC}})}.\hskip 1em
  plus 0.5em minus 0.4em\relax Maui, HI, USA: {IEEE}, Dec. 2012, pp.
  1325--1330.

\bibitem{Gomes2017c}
C.~Gomes, P.~Karalis, E.~M. {Navarro-L\'opez}, and H.~Vangheluwe,
  ``Approximated {{Stability Analysis}} of {{Bi}}-modal {{Hybrid
  Co}}-simulation {{Scenarios}},'' in \emph{1st {{Workshop}} on {{Formal
  Co}}-{{Simulation}} of {{Cyber}}-{{Physical Systems}}}.\hskip 1em plus 0.5em
  minus 0.4em\relax Trento, Italy: {Springer, Cham}, 2017, pp. 345--360.

\bibitem{Giannakopoulos2001}
F.~Giannakopoulos and K.~Pliete, ``Planar systems of piecewise linear
  differential equations with a line of discontinuity,'' \emph{Nonlinearity},
  vol.~14, no.~6, pp. 1611--1632, Oct. 2001.

\bibitem{Papachristodoulou2005}
A.~Papachristodoulou and S.~Prajna, ``A tutorial on sum of squares techniques
  for systems analysis,'' in \emph{American {{Control Conference}}}.\hskip 1em
  plus 0.5em minus 0.4em\relax Portland, OR, USA: {IEEE}, 2005, pp. 2686--2700.

\bibitem{proskurnikovdoes}
A.~V. Proskurnikov, ``Does sample-time emulation preserve exponential
  stability?'' 2020.

\bibitem{Arapostathis2007}
A.~Arapostathis and M.~E. Broucke, ``Stability and controllability of planar,
  conewise linear systems,'' vol.~56, no.~2, pp. 150--158.

\bibitem{johansson1999regularization}
K.~H. Johansson, M.~Egerstedt, J.~Lygeros, and S.~Sastry, ``On the
  regularization of zeno hybrid automata,'' \emph{Systems \& control letters},
  vol.~38, no.~3, pp. 141--150, 1999.

\bibitem{Johansson1999}
------, ``On the regularization of {{Zeno}} hybrid automata,'' \emph{Systems \&
  Control Letters}, vol.~38, no.~3, pp. 141--150, Oct. 1999.

\bibitem{Goebel2012}
\BIBentryALTinterwordspacing
R.~Goebel, R.~G. Sanfelice, and A.~R. Teel, \emph{Hybrid Dynamical Systems:
  Modeling, Stability, and Robustness}.\hskip 1em plus 0.5em minus 0.4em\relax
  Princeton University Press, 2012. [Online]. Available:
  \url{http://www.jstor.org/stable/j.ctt7s02z}
\BIBentrySTDinterwordspacing

\bibitem{blondel1999complexity}
V.~D. Blondel and J.~N. Tsitsiklis, ``Complexity of stability and
  controllability of elementary hybrid systems,'' \emph{Automatica}, vol.~35,
  no.~3, pp. 479--489, 1999.

\bibitem{Dai2012}
X.~Dai, ``A {{Gel}}'fand-type spectral radius formula and stability of linear
  constrained switching systems,'' \emph{Linear Algebra and its Applications},
  vol. 436, no.~5, pp. 1099--1113, Mar. 2012.

\bibitem{Philippe2016}
M.~Philippe, R.~Essick, G.~E. Dullerud, and R.~M. Jungers, ``Stability of
  discrete-time switching systems with constrained switching sequences,''
  \emph{Automatica}, vol.~72, pp. 242--250, Oct. 2016.

\bibitem{Lagarias1995}
J.~C. Lagarias and Y.~Wang, ``The finiteness conjecture for the generalized
  spectral radius of a set of matrices,'' \emph{Linear Algebra and its
  Applications}, vol. 214, no.~0, pp. 17--42, 1995.

\bibitem{Gomes2018c}
C.~Gomes, B.~Legat, R.~Jungers, and H.~Vangheluwe, ``Minimally {{Constrained
  Stable Switched Systems}} and {{Application}} to {{Co}}-simulation,'' in
  \emph{{{IEEE Conference}} on {{Decision}} and {{Control}}}, Miami Beach, FL,
  USA, 2018, pp. 5676--5681.

\bibitem{Hirsch2012}
M.~W. Hirsch, S.~Smale, and R.~L. Devaney, \emph{Differential Equations,
  Dynamical Systems, and an Introduction to Chaos}.\hskip 1em plus 0.5em minus
  0.4em\relax {Academic press}, 2012.

\bibitem{Munkres2000}
J.~R. Munkres, \emph{Topology}, 2nd~ed.\hskip 1em plus 0.5em minus 0.4em\relax
  {Prentice Hall}, 2000.

\bibitem{Liberzon2012}
D.~Liberzon, \emph{Switching in {{Systems}} and {{Control}}}.\hskip 1em plus
  0.5em minus 0.4em\relax {Springer Science \& Business Media}, 2012.

\bibitem{Robinson1998}
C.~Robinson, \emph{Dynamical Systems: Stability, Symbolic Dynamics, and
  Chaos}.\hskip 1em plus 0.5em minus 0.4em\relax {CRC press}, 1998.

\bibitem{Cellier2006}
F.~E. Cellier and E.~Kofman, \emph{Continuous {{System Simulation}}}.\hskip 1em
  plus 0.5em minus 0.4em\relax {Springer Science \& Business Media}, 2006.

\bibitem{rackauckas2017differentialequations}
C.~Rackauckas and Q.~Nie, ``Differentialequations.jl--a performant and
  feature-rich ecosystem for solving differential equations in julia,''
  \emph{Journal of Open Research Software}, vol.~5, no.~1, 2017.

\end{thebibliography}
